\documentclass{m2an}

\def \qed {\hfill \rule{1ex}{1ex}}
\def \U {{\bf u}}
\def \x {{\bf x}}

\def \H {{\bf H}}

\def \V {{\bf v}}
\def \bt {{\widetilde{\hbox{\rm \bf b}}}}
\def \W {{\bf w}}
\def \N {{\mathbf{n}}}

\def \Ut {{\tilde {\mathbf{u}}}}
\def \Deltat {{\widetilde \Delta}}
\def \F {{\bf f }}
\def \P {{\bf P}}

\def \L {{\bf L}}

\def \b {{\hbox{\rm  b}}}

\def \Deltat{{{\mathbf{\widetilde \Delta}}}}
\def \Nabla {{\mathbf{\nabla}}}
\def \Del {{\hbox{\boldmath $\delta$ \unboldmath \!\!}}}
\def \psig {{\hbox{\boldmath $\psi$ \unboldmath \!\!}}}

\newtheorem{lem}{Lemma}[section]
\newtheorem{theo}{Theorem}[section]
\newtheorem{prop}{Proposition}[section]

 \numberwithin{equation}{section}

\begin{document}
\title{Stability of a finite volume scheme\\ for the incompressible fluids}
%

\author{S\'ebastien Zimmermann}\address{17 rue Barr\`eme - 69006 LYON. e-mail: \texttt{Sebastien.Zimmermann@ec-lyon.fr}}

\date{Received: 13 august  2007}
\begin{abstract} We introduce a finite volume scheme for the
two-dimensional incompressible Navier-Stokes equations. We use a
triangular mesh. The unknowns for the velocity and pressure are
respectively piecewise constant and affine. We use a projection
method to deal with the incompressibility constraint. We show that
the differential operators in the Navier-Stokes equations and their
discrete counterparts share similar properties. In particular we
state an inf-sup (Babu\v{s}ka-Brezzi) condition. Using these
properties we infer the stability of the scheme.
 \end{abstract}
\begin{resume}
Nous introduisons ici un sch\'ema volumes finis pour les \'equations
de Navier-Stokes incompressibles en deux dimensions. Les maillages
consid\'er\'es sont form\'es de triangles. Les inconnues associ\'ees
\`a la vitesse et la pression sont respectivement constantes et
affines par morceaux. Nous utilisons une m\'ethode de projection
pour traiter la contrainte d'incompressibilit\'e. Nous v\'erifions
que les op\'erateurs diff\'erentiels apparaissant dans les
\'equations de Navier-Stokes et leurs analogues discrets v\'erifient
des propri\'et\'es similaires. Nous prouvons en particulier une
condition inf-sup. Nous en d\'eduisons la stabilit\'e du sch\'ema.
 \end{resume}
\subjclass{76D05, 74S10, 65M12}
\keywords{Incompressible fluids, Navier-Stokes equations, projection
methods, finite volume.}
\maketitle
\section{Introduction}

We consider the flow of an incompressible fluid in a 
polyhedral set  $\Omega \subset \mathbb{R}^2$ during the time
interval $[0,T]$. The velocity field $\mathbf{u}:\Omega \times [0,T]
\to \mathbb{R}^2 $ and the pressure field $p:\Omega \times [0,T] \to
\mathbb{R}$ satisfy the Navier-Stokes equations
\begin{eqnarray}
&& \mathbf{u}_t - \frac{1}{\hbox{Re}} \, \mathbf{\Delta} \mathbf{u}
+(\U \cdot \Nabla) \U+\nabla p = \F \, ,
\label{eq:mom} \\
&& \hbox{div } \U=0 \, , \label{eq:incomp}
\end{eqnarray}
with the boundary and initial conditions
\begin{equation*}
 \U|_{\partial \Omega}=0 \, , \hspace{2cm}  \U|_{t=0}=\U_0.
\end{equation*}
The terms $\mathbf{\Delta} \U$ and $(\U \cdot \Nabla)\U$ are
associated with the  physical phenomena of diffusion and convection,
respectively. The Reynolds number $\hbox{Re}$ measures the influence
of convection in the flow. For equations
(\ref{eq:mom})--(\ref{eq:incomp}), finite element and finite
difference methods are well known and mathematical studies are
available (see \cite{giraultr} for example). For finite volume
schemes, numerous computations have been conducted (\cite{kimchoi}
and \cite{boicaya} for example).
 However, few mathematical results are available in this case.  Let us cite {\sc Eymard and Herbin} \cite{herb3}
 and {\sc Eymard, Latch\'e and Herbin} \cite{eymard}.
 In order to deal with the incompressibility constraint (\ref{eq:incomp}), these works use a penalization  method.
 Another way is to use  the projection methods which have been
 introduced by {\sc Chorin} \cite{chorin} and {\sc Temam} \cite{temam}.
 This is the case in  {\sc Faure} \cite{faure} where
 the mesh is made of  squares.
 In {\sc Zimmermann} \cite{zimm1} the mesh is made of triangles, so that more
 complex geometries can be considered. 
 In the present paper the mesh is also made of triangles, but we consider a different
 discretization for the pressure. It leads to a linear system with a better-conditioned
 matrix.
The layout of the article is the following. We first introduce in
section \ref{sec:defd} the discrete setting. We state (section
\ref{subsec:defmaillage}) some notations and hypotheses on the mesh.
We define (section \ref{subsec:espd}) the spaces we use to
approximate the velocity and pressure. We define also (section
\ref{subsec:opd}) the operators we use to approximate the
differential operators in (\ref{eq:mom})--(\ref{eq:incomp}).
 Combining this with a projection method, we build the scheme in section
 \ref{sec:presschema}.
  In order to provide a mathematical analysis, we
  show in section \ref{sec:propop} that the differential operators
   in (\ref{eq:mom})--(\ref{eq:incomp}) and their discrete counterparts share similar properties.
 In particular, the discrete operators for the gradient and the
 divergence are adjoint. The discrete operator for the convection term
 is positive, stable and consistent. The discrete operator for the divergence
  satisfy an inf-sup (Babu\v{s}ka-Brezzi) condition.
From these properties we deduce in section \ref{sec:stab} the
stability of the scheme.

We conclude with  some notations. The spaces $(L^2,|.|)$ and
$(L^\infty,\|.\|_\infty)$ are the usual Lebesgue spaces and we set
$L^2_0=\{q \in L^2 \, ; \int_\Omega q(\x) \, d\x=0 \}$. Their
vectorial counterparts are $(\L^2,|.|)$ and
$(\L^\infty,\|.\|_\infty)$ with $\L^2=(L^2)^2$ and
$\L^\infty=(L^\infty)$. For $k \in \mathbb{N}^*$,
$(H^k,\|\cdot\|_k)$ is the usual Sobolev space. Its vectorial
counterpart is $(\H^k,\|.\|_k)$ with $\H^k=(H^k)^2$. For $k=1$, the
functions of $\H^1$ with a null trace on the boundary form the space
$\H^1_0$. Also, we set $\Nabla \U=(\nabla u_1,\nabla u_2)^T$ if
$\U=(u_1,u_2) \in \H^1$.
If $\mathbf{X}\subset \L^2$ is a Banach space, we define
${\mathcal{C}}(0,T;\mathbf{X})$
 (resp. $L^2(0,T;\mathbf{X})$) as the set of the
 applications $\mathbf{g}:[0,T] \to \mathbf{X}$ such that
$t \to |\mathbf{g}(t)|$ is continuous (resp. square integrable).
 The norm $\|.\|_{{\mathcal{C}}(0,T;\mathbf{X})}$
is defined by
$\|\mathbf{g}\|_{{\mathcal{C}}(0,T;\mathbf{X})}=\sup_{s\in[0,T]}
|\mathbf{g}(s)|$.
In all calculations, $C$ is a generic positive constant, depending
only on $\Omega$, $\U_0$ and $\F$.

\section{Discrete setting}
\label{sec:defd}

First, we introduce the spaces and the operators needed to build the
scheme.

\subsection{The mesh}
\label{subsec:defmaillage}
 Let ${\cal{T}}_h$ be a triangular mesh of $\Omega$.
The circumscribed circle of a triangle $K \in {\mathcal{T}}_h$ is
centered at  $\x_K$ and has the diameter $h_K$. We set $h=\max_{K
\in {\mathcal{T}}_h} h_K$. We assume that all the interior angles of
the triangles of the mesh are less than $\frac{\pi}{2}$, so that
$\x_K \in K$. The set of the edges of the triangle $K \in
{\mathcal{T}}_h$ is ${\cal{E}}_K$. The symbol $\N_{K,\sigma}$
denotes the unit vector normal to an edge $\sigma \in {\cal{E}}_K$
and  pointing outward $K$. We denote by ${\cal{E}}_h$ the set of the
edges of the mesh.
 We distinguish the subset ${\mathcal{E}}^{int}_h \subset {\mathcal{E}}_h$ (resp. ${\mathcal{E}}^{ext}_h$)
 of the edges located inside  $\Omega$ (resp. on $\partial \Omega$).
The middle of an  edge $\sigma \in {\cal{E}}_h$ is $\x_\sigma$ and
its length $|\sigma|$. For each edge $\sigma \in {\cal{E}}^{int}_h$,
 let $K_\sigma$ and $L_\sigma$ be the two  triangles having $\sigma$ in
 common. We set
 $d_\sigma=d(\x_{K_\sigma},\x_{L_\sigma})$. For all $\sigma \in {\mathcal{E}}^{ext}_h$, only the triangle $K_\sigma$
  located inside $\Omega$ is defined and we set $d_\sigma=d(\x_{K_\sigma},\x_{\sigma})$. Then for all $\sigma \in {\mathcal{E}}_h$ we set
   $\tau_\sigma=\frac{|\sigma|}{d_\sigma}$.
As in \cite{eymgal} we assume the following  on the mesh: there
exists $C>0$ such that
\begin{equation*}
   \forall \, \sigma \in {\cal{E}}_h \,,  \hspace{1cm}  d_\sigma \ge C \, |\sigma|
  \quad \mathrm{and} \quad
 |\sigma| \ge C \, h.
\end{equation*}
It implies that there exists $C>0$ such that
 \begin{equation}
 \label{eq:mintaus}
  \forall \, \sigma \in {\cal{E}}^{int}_h \, , \hspace{1cm} \tau_\sigma=|\sigma|/d_\sigma \ge C.
 \end{equation}

\subsection{The discrete spaces}

\label{subsec:espd}

\noindent   We first define
 \begin{equation*}
 \label{eq:defp0}
   P_0= \{ q \in L^2 \; ; \; \forall \, K \in {\mathcal{T}}_h, \; \; q|_K \hbox{ is a constant} \} \, ,
   \hspace{1cm} \P_0=(P_0)^2.
\end{equation*}
For the sake of concision, we set for all $q_h \in P_0$ (resp. $\V_h
\in \P_0$) and all triangle $K \in {\mathcal{T}}_h$:
 $q_K=q_h|_K$ (resp. $\V_K=\V_h|_K$).
Although $\P_0 \not \subset \H^1$, we define the discrete equivalent
of a $\H^1$ norm as follows. For all $\V_h \in \P_0$ we set
\begin{equation}
\label{eq:defh1d}
  \|\V_h\|_h=\left( \sum_{\sigma \in {\mathcal{E}}^{int}_h} \tau_{\sigma} \, |\V_{L_\sigma} - \V_{K_\sigma}|^2
  +\sum_{\sigma \in {\mathcal{E}}^{ext}_h} \tau_\sigma \, |\V_{K_\sigma}|^2
  \right)^{1/2}.
\end{equation}
 We have \cite{eymgal} a Poincar\'e-like inequality: there exists $C>0$ such that for all
  $\V_h \in \P_0$
 \begin{equation}
 \label{eq:inpoinp0}
  |\V_h| \le C \, \|\V_h\|_h.
 \end{equation}
We also have \cite{zimm1} an inverse inequality: there exists $C>0$
such that for all $\V_h \in \P_0$
\begin{equation}
\label{eq:ininvp0}
 h \, \|\V_h\|_h \le C \, |\V_h|.
\end{equation}
From the norm $\|.\|_h$ we deduce a dual norm.
 For all $\V_h \in \P_0$ we set
\begin{equation}
\label{eq:defnormdp0}
  \|\V_h\|_{-1,h}=\sup_{\psig_h \in \P_0} \frac{(\V_h,\psig_h)}{\|\psig_h\|_h}.
\end{equation}
For all $\U_h \in \P_0$ and $\V_h \in \P_0$ we have
  $(\U_h,\V_h) \le \|\U_h\|_{-1,h} \, \|\V_h\|_h$.
 We define the projection operator  $\Pi_{\P_0}: \L^2 \to \P_0$ as follows.
 For all $\W \in \L^2$, $\Pi_{\P_0} \W \in \P_0$ is given by
\begin{equation}
\label{eq:defprojp0}
   \forall \, K \in {\mathcal{T}}_h \, ,  \hspace{1cm} (\Pi_{\P_0} \W)|_K=\frac{1}{|K|} \int_K \W(\x) \, d\x.
\end{equation}
We easily check that for all  $\W \in \L^2$ and $\V_h \in \P_0$ we
have $(\Pi_{\P_0} \W, \V_h)=(\W,\V_h)$. We deduce from this that
$\Pi_{\P_0}$ is stable for the $\L^2$ norm.
 We define also the operator $\widetilde \Pi_{\P_0}:\H^2 \to \P_0$.
 For all $\W \in \H^2$, $\widetilde \Pi_{\P_0} \W \in \P_0$
 is given by
\begin{equation*}
  \forall \, K \in {\mathcal{T}}_h \, , \hspace{1cm} \widetilde \Pi_{\P_0} \W|_K=\W(\x_K).
\end{equation*}
According to the Sobolev embedding theorem, $\W \in \H^2$ is a.e.
equal to a continuous function. Therefore the definition above makes
sense. 
\noindent  We introduce also the finite element spaces
\begin{eqnarray}
\label{eq:defp1nc}
  P^d_1 &=& \{ v \in L^2 \; ; \; \forall \, K \in {\mathcal{T}}_h, \; \; v|_K \hbox{ is affine} \} \, , \nonumber\\
  P^{nc}_1 &=& \{ v_h \in P^d_1 \; ; \; \forall \, \sigma \in {\mathcal{E}}^{int}_h, \,
  v_h|_{K_\sigma}(\x_\sigma)=v_h|_{L_\sigma}(\x_\sigma) \, , \nonumber \\
  \P^{c}_1 &=& \{ \V_h \in (P^{d}_1)^2 \, ; \; \V_h \hbox{ is continuous} \hspace{.1cm} \hbox{ and } \hspace{.1cm}
  \V_h|_{\partial \Omega}=\mathbf{0}\}. \nonumber
 \end{eqnarray}
\noindent We have $\P^c_1 \subset \H^1_0$. We define
$\Pi_{\P^{c}_1}:\H^1_0 \to \P^{c}_1$. For all $\V=(v_1,v_2) \in
\H^1_0$,
 $\Pi_{\P^{c}_1} \V=(v^1_h,v^2_h) \in \P^{c}_1$ is given by
\begin{equation*}
  \forall \, \hbox{\boldmath $\phi$ \unboldmath}_{\hspace{-1.5mm} h}=(\phi^1_h,\phi^2_h)  \in \P^c_1 \, ,\hspace{1cm}   \sum_{i=1}^2 \big( \nabla v^i_h,\nabla \phi^i_h)
  = \sum_{i=1}^2 \big( \nabla v_i,\nabla \phi^i_h).
\end{equation*}
The operator $\Pi_{\P^c_1}$
 is stable for the  $\H^1$ norm. One checks (\!  \cite{brenns} p. 110) that there exists    $C>0$ such that for
all $\V \in \H^1$
\begin{equation}
\label{eq:errintp1c}
   |\V-\Pi_{\P^c_1} \V| \le C \, h \, \|\V\|_1.
\end{equation}
Let us address now the space  $P^{nc}_1$. If $q_h \in P^{nc}_1$, we
have usually $\nabla q_h \not \in \L^2$. Thus we define the operator
$\nabla_h: P^{nc}_1 \hspace{-.2cm}\to \P_0 $
 by setting  for all $q_h \in P_0$ and all triangle  $K \in {\mathcal{T}}_h$
\begin{equation}
\label{eq:defg}
  \nabla_h q_h |_K = \frac{1}{|K|} \, \int_K \nabla q_h \, d\x.
\end{equation}
The associated norm is defined by
\begin{equation*}
\label{eq:defnormp1nc}
  \|q_h\|_{1,h}=\left(|q_h|^2+|\nabla_h q_h|^2 \right)^{1/2}.
\end{equation*}
We have a  Poincar\'e-like inequality : there exists $C>0$ such that
for all $q_h \in P^{nc}_1 \cap L^2_0$
\begin{equation}
\label{eq:inpoinp1nc}
  |q_h| \le C \, |\nabla_h q_h|.
\end{equation}
We define the projection operator $\Pi_{P^{nc}_1}$. For all $q \in
H^1$, $\Pi_{P^{nc}_1} q$ is given by
\begin{equation*}
\label{eq:defpp1nc}
   \forall \, \sigma \in {\mathcal{E}}_h \, , \hspace{1cm}
    \int_\sigma (\Pi_{P^{nc}_1} q) \, d\sigma=\int_\sigma q \, d\sigma.
\end{equation*}
One checks (\! \cite{brenns} p.110) that there exists $C>0$ such
that
\begin{equation}
\label{eq:estintp1nc}
  |p-\Pi_{P^{nc}_1}p| \le C \, h\, \|p\|_1 \, , \hspace{1cm} \left|\widetilde \nabla_h(p-\Pi_{P^{nc}_1}p)\right| \le C \,
  \|p\|_1.
\end{equation}

\noindent  Finally, we  use the Raviart-Thomas spaces (see
\cite{brezzfor})
\begin{eqnarray*}
\label{eq:defrt0}
  \mathbf{RT^d_0}&=& \{ \V_h \in \P^d_1 \; ; \quad \forall \, \sigma \in {\mathcal{E}}_K,
   \quad \V_h|_K \cdot \N_{K,\sigma} \hbox{ is a  constant,}  \hspace{.3cm} \hbox{ and }  \hspace{.3cm}
    \V_h \cdot \N|_{\partial \Omega}=0  \} \, , \nonumber \\
  \mathbf{RT_0}&=& \{ \V_h \in \mathbf{RT^d_0} \; ; \quad \forall \, K \in {\mathcal{T}}_h,
   \quad \forall \, \sigma \in {\mathcal{E}}_K,
 \quad   \V_h|_{K_\sigma} \cdot \N_{K_\sigma,\sigma} = \V_h|_{L_\sigma} \cdot \N_{K_\sigma,\sigma}\}.
\end{eqnarray*}
For all $\V_h \in \mathbf{RT_0}$, $K \in {\mathcal{T}}_h$ and
$\sigma \in {\mathcal{E}}_K$ we set
  $(\V_h \cdot \N_{K,\sigma})_\sigma=\V_h|_K \cdot \N_{K,\sigma}$.
\noindent We define the operator  $\Pi_{\mathbf{RT_0}}:\H^1 \to
\mathbf{RT_0}$. For all $\V \in \H^1$, $\Pi_{\mathbf{RT_0}} \V \in
\mathbf{RT_0}$ is given by
\begin{equation}
\label{eq:defprt0} \forall \, K \in {\mathcal{T}}_h \, ,
\hspace{.5cm}
  \forall \, \sigma\in {\mathcal{E}}_K \, , \hspace{1cm}
    (\Pi_{\mathbf{RT_0}} \V \cdot \N_{K,\sigma})_\sigma=\frac{1}{|\sigma|}\int_\sigma \V \, d\sigma.
\end{equation}

\subsection{The discrete operators}
\label{subsec:opd}

The equations (\ref{eq:mom})--(\ref{eq:incomp}) use the differential
operators gradient, divergence and laplacian.
 Using the spaces of section \ref{subsec:espd}, we define their discrete counterparts.
The discrete gradient $ \nabla_h: P^{nc}_1 \to \P_0$ is defined by
(\ref{eq:defg}). \noindent  The discrete divergence operator
$\hbox{div}_h: \P_0 \to P^{nc}_1$ is built so that it is adjoint to
the operator $\nabla_h$. We set for all $\V_h \in \P_0$ and all
triangle $K \in {\mathcal{T}}_h$
 \begin{eqnarray}
\label{eq:defdivh1} &&  \forall \, \sigma \in {\mathcal{E}}^{int}_h,
\hspace{2cm}  (\hbox{div}_h \, \V_h) (\x_\sigma)=
\frac{3 \, |\sigma|}{|K_\sigma|+|L_\sigma|} \, (\V_{L_\sigma} - \V_{K_\sigma}) \cdot \N_{K,\sigma} \, ; \nonumber \\
&&  \forall \, \sigma \in {\mathcal{E}}^{ext}_h, \hspace{2cm}
(\hbox{div}_h \, \V_h) (\x_\sigma)= -\frac{3 \,
|\sigma|}{|K_\sigma|+|L_\sigma|} \, \V_{K_\sigma} \cdot
\N_{K,\sigma}.  \label{eq:defdivh2}
\end{eqnarray}
The first discrete laplacian $\Delta_h:P^{nc}_1 \to P^{nc}_1$
ensures that the incompressibility constraint (\ref{eq:incomp}) is
satisfied in a discrete sense (see the proof of proposition
\ref{prop:umrt0} below). We set for all $q_h \in P^{nc}_1$
\begin{equation*}
\label{eq:deflap} \Delta_h q_h=\hbox{div}_h (\nabla_h q_h).
\end{equation*}
\noindent The second discrete laplacian
$ \Deltat_h: \P_0 \to \P_0$ is the usual operator in finite volume
schemes \cite{eymgal}. We set for all $\V_h \in \P_0$ and all
triangle $K \in {\mathcal{T}}_h$
\begin{equation*}
\label{eq:defl}
  \Deltat_h \V_h|_K
  = \frac{1}{|K|} \sum_{\sigma \in {\mathcal{E}}_K \cap {\mathcal{E}}^{int}_h}
 \tau_\sigma \, (\V_{L_\sigma} - \V_{K_\sigma} )
-\frac{1}{|K|}\sum_{\sigma \in {\mathcal{E}}_K \cap
{\mathcal{E}}^{ext}_h}
 \tau_\sigma \, \V_{K_\sigma}.
\end{equation*}


\vspace{.1cm}

\noindent In order to approximate the term  $(\U \cdot \Nabla) \U$
in (\ref{eq:mom}) we define a bilinear form $\bt_h: \mathbf{RT_0}
\times \P_0 \to \P_0$ using the well-known upwind scheme
\cite{eymgal}. For all $\U_h \in \P_0$, $\V_h \in \P_0$, and all
triangle $K \in {\mathcal{T}}_h$ we set
\begin{equation}
\label{eq:defbth}
  \bt_h(\U_h,\V_h)\big|_K=
\frac{1}{|K|} \sum_{\sigma \in {\mathcal{E}}_K \cap
{\mathcal{E}}^{int}_h} \,
   |\sigma| \, \Big( (\U \cdot \N_{K,\sigma})^+_\sigma \, \V_K +
(\U \cdot \N_{K,\sigma})^-_\sigma \, \V_{L_\sigma}  \Big).
\end{equation}
We have set $a^+=\max(a,0)$, $a^-=\min(a,0)$ for all
$a\in\mathbb{R}$. Lastly, we define the trilinear form $ \b_h:
\mathbf{RT_0} \times \P_0 \times \P_0 \to \mathbb{R}^2$ as follows.
For all $\U_h \in \mathbf{RT_0}$, $\V_h \in \P_0$, $\W_h \in \P_0$,
we set
\begin{equation}
\label{eq:defbh}
  \b_h(\U_h,\V_h,\W_h)=\sum_{K \in {\mathcal{T}}_h} |K| \, \W_K \cdot
   \bt_h(\U_h,\V_h)\big|_K .
\end{equation}

\section{The scheme}

\label{sec:presschema}

In order to deal with the incompressibility constraint
(\ref{eq:incomp}) we use a projection method. This kind of method
has been introduced by {\sc Chorin} \cite{chorin}
 and {\sc Temam} \cite{temam}. The basic idea is the following.
The time interval $[0,T]$ is split with a time step $k$:
$[0,T]=\bigcup_{n=0}^N [t_n,t_{n+1}]$ with $N \in \mathbb{N}^*$ and
$t_n=n \, k$ for all $n\in\{0,\dots,N\}$. For all
$m\in\{2,\dots,N\}$, we compute  (see equation (\ref{eq:mombdf})
below) a first velocity field $\Ut^m_h \simeq \U(t_m)$ using only
equation (\ref{eq:mom}). We use a second-order BDF scheme for the
discretization in time. We then project $\Ut^m_h$ (see equation
(\ref{eq:projib}) below) over a subspace of $\P_0$. We get a a
pressure field  $p^m_h \simeq p(t_m)$ and a second velocity field
$\U^m_h \simeq \U(t_m)$, which fulfills the incompressibility
constraint (\ref{eq:incomp}) in a discrete sense. The algorithm goes
as follows. For all $m\in\{0,\dots,N\}$, we set $\F^m_h=\Pi_{\P_0}
\F(t_m)$. Since the operator $\Pi_{\P_0}$ is stable for the
$\L^2$-norm we get
\begin{equation}
\label{eq:propfh}
  |\F^m_h|=|\Pi_{\P_0} \F(t_m)| \le |\F(t_m)| \le \|\F\|_{{\mathcal{C}}(0,T;\L^2)}.
\end{equation}
We start with the initial values
\begin{equation*}
  \U^0_h \in \P_0 \cap \mathbf{RT_0} \, , \hspace{1cm}   \U^1_h \in \P_0 \cap \mathbf{RT_0}
  \,\hspace{1cm}   p^1_h \in P_0 \cap L^2_0.
\end{equation*}
For all $n \in \{1,\dots,N\}$, $(\Ut^{n+1}_h,p^{n+1}_h,\U^{n+1}_h)$
is deduced from $(\Ut^n_h,p^n_h,\U^n_h)$ as follows.

\begin{itemize}
\item $\Ut^{n+1}_h \in \P_0$ is given by
\begin{equation}
\label{eq:mombdf} \hspace{-.2cm} \frac{3 \, \Ut^{n+1}_h -4 \,
\U^n_h+\U^{n-1}_h}{2 \, k} - \frac{1}{\hbox{Re}} \, \Deltat_h
\Ut^{n+1}_h+\bt_h(2\, \U^n_h-\U^{n-1}_h, \Ut^{n+1}_h)+\nabla_h p^n_h
=\F^{n+1}_h  \, ,
\end{equation}

\item $p^{n+1}_h \in P^{nc}_1 \cap L^2_0$ is the solution of
\begin{equation}
\label{eq:projia} \Delta_h (p^{n+1}_h-p^n_h)=\frac{3}{2 \, k}\,
\hbox{div}_h \, \Ut^{n+1}_h  \, ,
\end{equation}

\item $\U^{n+1}_h \in  \P_0$ is deduced by
\begin{equation}
\label{eq:projib}
 \U^{n+1}_h = \Ut^{n+1}_h -\frac{2\, k}{3} \, \nabla_h (p^{n+1}_h - p^n_h).
\end{equation}
\end{itemize}
Existence and unicity of a solution to equation (\ref{eq:mombdf})
   is classical (\cite{eymgal} for example).
The convection term in (\ref{eq:mombdf}) is well defined thanks to
the following result.


\begin{prop}
\label{prop:umrt0} For all $m\in\{0,\dots,N\}$ we have  $\U^m_h
\in\mathbf{RT_0}$ .
\end{prop}
\noindent {\sc Proof.} If $m\in\{0,1\}$  the result holds by
definition.
 If $m\in\{2,\dots,N\}$ we apply the operator $\hbox{div}_h$ to
(\ref{eq:projia}) and compare with (\ref{eq:projib}). We get
$\hbox{div}_h \, \U^m_h=0$. Using definition (\ref{eq:defdivh1}) we
get $\U^m_h \in \mathbf{RT_0}$. \qed

\vspace{.2cm}

\noindent  Let us show
  that equation (\ref{eq:projia}) also has  a unique solution.
Let  $q_h \in P^{nc}_1 \cap L^2_0$ such that $\Delta_h q_h=0$.
 According to proposition \ref{prop:propadjh} we have for all $q_h \in P_0$
\begin{equation*}
  -(\Delta_h q_h,q_h)=-\big(\hbox{div}_h(\nabla_h q_h),q_h\big)=(\nabla_h q_h,\nabla_h q_h)=|\nabla_h q_h|^2.
\end{equation*}
Therefore we have  $\nabla_h q_h=0$, so that $q_h=0$ since $q_h \in
L^2_0$. We have thus proved the unicity of a solution for
(\ref{eq:projia}). It is also the case for the associated linear
system.
 It implies that this linear system has indeed a solution. Hence it is also the case for
 equation  (\ref{eq:projia}).
Note finally that since $\U^m_h \in \P_0 \cap \mathbf{RT_0}$, we
have $\hbox{div} \, \U^m_h=0$ for all $m\in\{0,\dots,N\}$. Hence the
incompressibility condition (\ref{eq:incomp}) is fulfilled.


\section{Properties of the discrete operators}
\label{sec:propop}

We show that the differential operators in
(\ref{eq:mom})--(\ref{eq:incomp}) and the operators defined in
section \ref{subsec:opd} share similar properties.

\subsection{Properties of the discrete convective term}
\label{subsec:propconvh}

We define   $\bt: \H^1 \times \H^1 \to \L^2$. For all
 $\U \in \H^1$ and $\V=(v_1,v_2) \in \H^1$ we set
$ \label{eq:defbt} \bt(\U,\V)= \big(\hbox{div}(v_1 \,
\U),\hbox{div}(v_2 \, \U)\big)$. We show that the operator $\bt_h$
is a consistent approximation of $\bt$.
\begin{prop}
\label{prop:consbh} There exists a constant $C>0$ such that for all
$\V \in \H^{2}$ and all $\U \in \H^2 \cap \H^1_0$ satisfying
$\hbox{\rm div} \, \U=0$
\begin{equation*}
\label{eq:consbh}
  \|\Pi_{\P_0}\bt(\U,\V) - \bt_h(\Pi_{\mathbf{RT_0}} \U, \widetilde \Pi_{\P_0} \V)\|_{-1,h}
\le C \, h  \,  \|\U\|_2 \,  \|\V\|_{1}.
\end{equation*}
\end{prop}
\noindent {\sc Proof.}  We set  $\U_h=\Pi_{\mathbf{RT_0}} \U$ and
$\V_h=\widetilde \Pi_{\P_0} \V$.
%
Let $K \in {\mathcal{T}}_h$. According to the  divergence formula
and (\ref{eq:defprojp0}) we have
\begin{equation*}
\label{eq:eqconsb2}
 \Pi_{\P_0}\bt(\U,\V) |_K
 =\frac{1}{|K|} \sum_{\sigma \in {\mathcal{E}}_K \cap {\mathcal{E}}^{int}_h} \int_\sigma \V \, (\U \cdot \N) \, d\sigma.
\end{equation*}
On the other hand, let us rewrite $\bt_h(\U_h,\V_h)$. Let $\sigma
\in {\mathcal{E}}_K \cap {\mathcal{E}}^{int}_h$. Setting
\begin{equation*}
  \V_{K,L_\sigma}= \left\{
\begin{array}{ll}
   \V_{K} &\hbox{ si } (\U_h \cdot \N_{K,\sigma})_\sigma \ge 0 \\
   \V_{L_\sigma} &\hbox{ si } (\U_h \cdot \N_{K,\sigma})_\sigma < 0
\end{array}
\right.
\end{equation*}
one  checks that
  $\V_K \, (\U_h \cdot \N_{K,\sigma})^+_\sigma + \V_{L_\sigma} \,  (\U_h \cdot
  \N_{K,\sigma})^-_\sigma
=  \V_{K,L_\sigma} \, (\U_h \cdot \N_{K,\sigma})_\sigma$.
 Using (\ref{eq:defprt0}),  we deduce from
(\ref{eq:defbth}) that
\begin{equation*}
\label{eq:eqconsb3} \bt_h(\U_h,\V_h)|_K=\frac{1}{|K|} \sum_{\sigma
\in {\mathcal{E}}_K \cap {\mathcal{E}}^{int}_h} \int_\sigma
\V_{K,L_\sigma} \, (\U_h \cdot \N_{K,\sigma}) \, d\sigma.
\end{equation*}
Thus
\begin{equation*}
\label{eq:eqconsb4}
  \big( \Pi_{\P_0}\bt(\U,\V)-\bt_h(\U_h,\V_h) \big)|_K
 =\frac{1}{|K|} \sum_{\sigma \in {\mathcal{E}}_K \cap {\mathcal{E}}^{int}_h}
 \int_\sigma (\V-\V_{K,L_\sigma}) \, (\U_h \cdot \N) \, d\sigma.
\end{equation*}
 Let $\psig_h \in \P_0$. We have
\begin{eqnarray*}
\label{eq:eqconsb5}
 \hspace{-.1cm} \Big( \Pi_{\P_0}\bt(\U,\V) - \bt_h(\U_h, \V_h), \psig_h \Big)
&=&\sum_{K \in {\mathcal{T}}_h}   \psig_K  \sum_{\sigma \in
{\mathcal{E}}_K \cap {\mathcal{E}}^{int}_h}
 \int_\sigma (\V-\V_{K,L_\sigma}) \, (\U_h \cdot \N) \, d\sigma \nonumber \\
 &=& \sum_{\sigma \in {\mathcal{E}}^{int}_h} (\psig_{K_\sigma}-\psig_{L_\sigma})  \int_\sigma (\V-\V_{K_\sigma,L_\sigma}) \, (\U_h \cdot \N) \, d\sigma.
\end{eqnarray*}
Let $\sigma \in {\mathcal{E}}^{int}_h$.  We consider the
quadrilateral $D_\sigma$ defined by $\x_{K_\sigma}$, $\x_{L_\sigma}$
and  the vertex of $\sigma$.
 We set
\begin{equation*}
  D_{K,L_\sigma}= \left\{
\begin{array}{ll}
   D_{\sigma} \cap K &\hbox{ si } (\U_h \cdot \N_{K,\sigma})_\sigma \ge 0 \\
   D_{\sigma} \cap L_\sigma &\hbox{ si } (\U_h \cdot \N_{K,\sigma})_\sigma < 0
\end{array}
\right. .
\end{equation*}
Using a Taylor expansion and a density argument (see \cite{zimm1})
one checks that
\begin{equation*}
\label{eq:eqconsb7c} \int_\sigma |\V-\V_{K_\sigma,L_\sigma}| \,
d\sigma \le
 C \, h\left(  \int_{D_{K_\sigma,L_\sigma}} |\Nabla \V \, (\mathbf{y})|^2  \,  d\mathbf{y} \right)^{1/2}.
\end{equation*}
Thus
\begin{eqnarray*}
&& \left|  \big( \Pi_{\P_0}\bt(\U,\V) - \bt_h(\Pi_{\mathbf{RT_0}}
\U, \widetilde \Pi_{\P_0} \V), \psig_h \big) \right|
  \\
&& \le C \, h \, \|\U\|_{\H^{2}} \left( \sum_{\sigma \in
{\mathcal{E}}^{int}_h} |\psig_{L_\sigma} - \psig_{K_\sigma}|^2
\right)^{1/2}
 \left( \sum_{\sigma \in {\mathcal{E}}^{int}_h}   \int_{D_{K_\sigma,L_\sigma}} |\Nabla \V \, (\mathbf{y})|^2
  \,  d\mathbf{y} \right)^{1/2}
\end{eqnarray*}
so that
$\left|  \big( \Pi_{\P_0}\bt(\U,\V) - \bt_h(\Pi_{\mathbf{RT_0}} \U,
\widetilde \Pi_{\P_0} \V), \psig_h \big) \right|
  \le C \, h \,  \|\U\|_{\H^{2}} \, \|\psig_h\|_{1,h} \, \|\V\|_1$.
Using then definition (\ref{eq:defnormdp0}), we get the result. \qed

\vspace{.2cm}

\noindent Let  $\V \in \L^\infty \cap \H^1$ and $\U \in \H^1$ with
$\hbox{div} \, \U \ge 0$ a.e. in $\Omega$.
 Integrating by parts one checks that 
    $\int_\Omega \V \cdot \bt(\U,\V) \, d\x =\int_\Omega  \frac{|\V|^2}{2}\,  \hbox{div} \, \U\, d\x \ge 0$.
The  operator $\b_h$ shares a similar property.

\vspace{-.2cm}

\begin{prop}
\label{prop:posbh} Let $\U_h \in \mathbf{RT_0}$ such that $\hbox{\rm
div} \, \U_h \ge  0$. For all $\V_h \in \P_0$  we have
\begin{equation*}
\label{eq:bhpos}
  \b_h(\U_h,\V_h,\V_h) \ge 0.
\end{equation*}
\end{prop}
\noindent {\sc Proof.} Remember that for all edges $\sigma \in
{\mathcal{E}}^{int}_h$, two triangles $K_\sigma$ et $L_\sigma$ share
$\sigma$ as an edge. We denote by $K_\sigma$ the one such that
$\U_{\sigma} \cdot \N_{K_\sigma,\sigma} \ge 0$. Using the algebraic
identity $2 \, a \, (a-b)=a^2-b^2+(a-b)^2$ we deduce from
(\ref{eq:defbh})
\begin{eqnarray*}
  2 \, \b_h(\U_h,\V_h,\V_h)&=&
2\sum_{\sigma \in {\mathcal{E}}^{int}_h} |\sigma| \, \V_{K\sigma}
\cdot (\V_{K_\sigma} - \V_{L_\sigma})
 \, (\U_{h} \cdot \N_{K_\sigma,\sigma}) \\
&=&\sum_{\sigma \in {\mathcal{E}}^{int}_h} |\sigma| \,  \Big(
|\V_{K\sigma}|^2 - |\V_{L_\sigma}|^2+|\V_{K_\sigma} -
\V_{L_\sigma}|^2 \Big) \, (\U_{h} \cdot \N_{K_\sigma,\sigma})
\end{eqnarray*}
so that
$2 \, \b_h(\U_h,\V_h,\V_h) \ge  \sum_{\sigma \in
{\mathcal{E}}^{int}_h} |\sigma|  \, \Big( |\V_{K\sigma}|^2 -
|\V_{L_\sigma}|^2 \Big) \, (\U_{h} \cdot \N_{K_\sigma,\sigma})$.
This sum can be written  as a sum over the triangles of the mesh. We
get
\begin{equation*}
2 \, \b_h(\U_h,\V_h,\V_h) \ge  \sum_{K \in {\mathcal{T}}_h}
|\V_{K_\sigma}|^2 \sum_{\sigma \in {\mathcal{E}}_K \cap
{\mathcal{E}}^{int}_h} |\sigma| \, (\U_h \cdot
\N_{K_\sigma,\sigma}).
\end{equation*}
Using finally the divergence formula we get
\begin{equation*}
\hspace{2.8cm} 2 \, \b_h(\U_h,\V_h,\V_h) \ge \sum_{K \in
{\mathcal{T}}_h} |K| \, |\V_{K}|^2 \,    \int_K \hbox{\rm div} \,
\U_h \, d\x \ge 0. \hspace{2.8cm} \qed
\end{equation*}

\vspace{-.2cm}

\noindent The following result states that the operator $\b_h$ is
stable for suitable norms.
\begin{prop}
\label{prop:stabbth} There exists a constant $C>0$ such that for all
 $\V_h \in \P_0$, $\W_h\in \P_0$, $\U_h \in \P_0$ satisfying  $\hbox{ \rm div} \, \U_h=0$
\begin{equation*}
\label{eq:majbh} |\b_h(\U_h,\V_h,\V_h)| \le C \, |\U_h| \,
\|\V_h\|_{h} \, \|\V_h\|_h.
\end{equation*}
\end{prop}
\noindent {\sc Proof.} For all triangle $K \in {\mathcal{T}}_h$ and
all edge  $\sigma \in {\mathcal{E}}_K \cap {\mathcal{E}}^{int}_h$,
we have
\begin{equation*}
  (\U_h \cdot \N_{K,\sigma})^+_\sigma \, \V_K +  (\U_h \cdot
  \N_{K,\sigma})^-_\sigma
  \, \V_{L_\sigma}
  =(\U_h \cdot \N_{K,\sigma})_\sigma \, \V_K
 - |(\U_h \cdot \N_{K,\sigma})_\sigma| \, (\V_{L_\sigma}-\V_K).
\end{equation*}
Using this splitting, we deduce from  (\ref{eq:defbh})
 $\b_h(\U_h,\V_h,\W_h) = S_1 +S_2$
with
\begin{eqnarray*}
 S_1 &=&
\sum_{K \in {\mathcal{T}}_h} \V_K \cdot \W_K \sum_{\sigma \in
{\mathcal{E}}_K \cap {\mathcal{E}}^{int}_h} |\sigma|
\,  (\U_h \cdot \N_{K,\sigma})_\sigma \, , \nonumber \\
S_2 &=& -\sum_{K \in {\mathcal{T}}_h} \W_K \cdot \sum_{\sigma \in
{\mathcal{E}}_K \cap {\mathcal{E}}^{int}_h} |\sigma|\, |(\U_h \cdot
\N_{K,\sigma})_\sigma| \, (\V_{L_\sigma}-\V_K) \label{eq:defs2}.
\end{eqnarray*}
By writing the sum over the edges as a sum over the triangles we
have
\begin{equation*}
S_2 = -\sum_{\sigma \in  {\mathcal{E}}^{int}_h}   |\sigma|\, |(\U_h
\cdot \N_{K,\sigma})_\sigma| \, (\V_{L_\sigma}-\V_K) \cdot
(\W_{L_\sigma}-\W_K).
\end{equation*}
Using the Cauchy-Schwarz inequality we get
\begin{equation*}
  |S_2| \le  h \, \|\U_h\|_{\infty} \, \left( \sum_{\sigma \in {\mathcal{E}}^{int}_h}
  |\V_{L_\sigma} - \V_{K_\sigma}|^2 \right)^{1/2} \,  \left( \sum_{\sigma \in {\mathcal{E}}^{int}_h} |\W_{L_\sigma} - \W_{K_\sigma}|^2 \right)^{1/2}.
\end{equation*}
Since $\U_h \in \mathbf{RT_0}$ we have \cite{eymgal} the inverse
inequality
  $h \, \|\U_h\|_{\infty} \le C \, |\U_h|$.
 Using (\ref{eq:mintaus}) and  (\ref{eq:defh1d}) we get
\begin{equation*}
 \sum_{\sigma \in {\mathcal{E}}^{int}_h} |\V_{L_\sigma} - \V_{K_\sigma}|^2
 \le  C \, \sum_{\sigma \in {\mathcal{E}}^{int}_h} \tau_\sigma \, |\V_{L_\sigma} - \V_{K_\sigma}|^2
 \le C \, \|\V_h\|^2_h
\end{equation*}
and in a similar way
 $\sum_{\sigma \in {\mathcal{E}}^{int}_h} |\W_{L_\sigma} - \W_{K_\sigma}|^2  \le C \, \|\W_h\|^2_h$.
Thus
  $|S_2| \le C \, |\U_h| \, \|\V_h\|_h \, \|\W_h\|_h$.
On the other hand, according to the divergence formula
\begin{equation*}
\label{eq:estbh1} S_1 = \sum_{K \in {\mathcal{T}}_h}  |K| \,  (\V_K
\cdot \W_K) \, \int_K \hbox{div} \, \U_h \, d\x=0.
\end{equation*}
By gathering the estimates for $S_1$ and $S_2$ we get the
result.\qed

\subsection{Properties of the discrete divergence}
\label{subsec:propdivh}

The operators gradient and divergence are adjoint: if $q \in H^1$ ,
$\V \in \H^1$ with $\V \cdot \N|_{\partial \Omega}=0$, we get
$(\V,\nabla q)=-(q,\hbox{div} \, \V)$ by integrating by parts. For
$\nabla_h$ and $\hbox{div}_h$ we state the following.

\vspace{-.2cm}

\begin{prop}
\label{prop:propadjh} For all $\V_h \in \P_0$ and $q_h \in P^{nc}_1$
we have: $\label{eq:opdadj}
  (\V_h, \nabla_h q_h)=-(q_h,\hbox{\rm div}_h \, \V_h)$.
\end{prop}
\noindent {\sc Proof.} According to (\ref{eq:defg})
\begin{equation*}
  (\V_h, \nabla_h q_h)   =
\sum_{K \in {\cal{T}}_h} |K| \, \V_K \cdot \nabla_h q_h|_K
 = \sum_{K \in {\cal{T}}_h} \V_K \cdot \Big( \sum_{\sigma \in {\cal{E}}_K} |\sigma|
q_h(\x_\sigma) \,\N_{K,\sigma} \Big).
\end{equation*}
By writing this sum as a sum over the edges we get
\begin{equation}
\label{eq:eqadj4} (\V_h, \nabla_h q_h)
  =  -\sum_{\sigma \in {\cal{E}}^{int}_h} |\sigma| \, q_h(\x_\sigma) \, (\V_{L_\sigma}-\V_{K_\sigma})
 \cdot \N_{{K_\sigma},\sigma}
+\sum_{\sigma \in {\cal{E}}^{ext}_h} |\sigma| \, q_h(\x_\sigma) \,
\V_{K_\sigma} \cdot \N_{K_\sigma,\sigma}.
\end{equation}
On the other hand, using a quadrature formula
\begin{equation*}
 -(q_h,\hbox{div}_h \, \V_h)
= -\sum_{K \in {\cal{T}}_h} \frac{|K|}{3} \sum_{\sigma \in
{\cal{E}}_K} q_h(\x_\sigma) \, (\hbox{div}_h \, \V_h) (\x_\sigma).
\end{equation*}
By writing this sum as a sum over the edges of the mesh we get
\begin{equation*}
-(q_h,\hbox{div}_h \, \V_h)=-\sum_{\sigma \in {\cal{E}}^{int}_h}
\Big( \frac{|K_\sigma|}{3} +\frac{|L_\sigma|}{3}\Big) \,
q_h(\x_\sigma) \, (\hbox{div}_h \V_h) (\x_\sigma)  -\sum_{\sigma \in
{\cal{E}}^{ext}_h} \frac{|K_\sigma|}{3} \, q_h(\x_\sigma) \,
(\hbox{div}_h \, \V_h) (\x_\sigma).
\end{equation*}
Using definition (\ref{eq:defdivh1}) and comparing with
(\ref{eq:eqadj4}) we get the result. \qed

\vspace{.2cm}

\noindent The divergence operator and the spaces $L^2_0$, $\H^1_0$
satisfy
 the following property, called inf-sup (or Babu\v{s}ka-Brezzi) condition (see \cite{giraultr} for example).
There exists a constant $C>0$ such that
\begin{equation}
\label{eq:infsc}
  \inf_{\begin{array}{l}
  {\scriptstyle q \in L^2_0 \backslash \{0\} }
  \end{array}}
\hspace{.2cm} \sup_{\V \in \H^1_0 \backslash \{\mathbf{0}\}}
-\frac{(q,\hbox{\rm div} \, \V)}{\|\V\|_1 |q|} \ge C.
\end{equation}
We will now show that the operator $\hbox{\rm div}_h$ and the spaces
$P_0 \cap L^2_0$, $\P_0$ satisfy an analogous property. The proof
uses the following lemma.


\begin{lem}
\label{lemma:infsd}
  There exists a
  constant $C>0$ such that
  \begin{equation*}\label{eq:infsdeg}
    \forall \, q_h \in P^{nc}_1 \cap L^2_0 \, , \hspace{1cm} \sup_{\V_h \in \P_0\backslash \{\mathbf{0}\}}
    -\frac{(q_h,\hbox{\rm div}_h \, \V_h)}{\|\V_h\|_h} \ge C \, h \, \|q_h\|_{1,h}.
  \end{equation*}
\end{lem}
\noindent {\sc Proof.} If $q_h=0$ the result is trivial. Let $q_h
\in P^{nc}_1 \cap L^2_0\backslash \{0\}$. Let $\V_h=\nabla_h q_h \in
\P_0\backslash \{ \mathbf{0}\}$. Using proposition
\ref{prop:propadjh} we have
\begin{equation*}
  -(q_h,\hbox{div}_h \V_h)=(\V_h,\nabla_h q_h)=|\nabla_h q_h|^2=|\nabla_h q_h| \, |\V_h|.
\end{equation*}
Using (\ref{eq:inpoinp0}) and (\ref{eq:ininvp0})  we get
$-(q_h,\hbox{div}_h \V_h) \ge C \, h  \, \|q_h\|_{1,h} \,
\|\V_h\|_h$.  \qed


\noindent We now state the result.


\begin{prop}
\label{prop:condinfs} There exists a constant $C>0$ such that for
all $q_h \in P^{nc}_1 \cap L^2_0$
\begin{equation*}
\label{eq:infsd}
   \sup_{\V_h \in \P_0 \backslash \{\mathbf{0}\}} -\frac{(q_h,\hbox{\rm div}_h \, \V_h)}{\|\V_h\|_h} \ge C
    \, |q_h|.
\end{equation*}
\end{prop}
\noindent {\sc Proof.} If $q_h=0$ the result is trivial. Let $q_h
\in P^{nc}_1 \cap L^2_0 \backslash \{0\}$. According to
(\ref{eq:infsc})
 there exists $\V \in \H^1_0$ such that
\begin{equation}
\label{eq:propv}
  \hbox{div} \, \V=-q_h \hspace{.3cm} \hbox{ and } \hspace{.3cm} \|\V\|_1
   \le C \, |q_h|.
\end{equation}
We set
  $\V_h=\Pi_{\P^c_1} \V$.
We want to  estimate $-\big(q_h,\hbox{div}_h(\Pi_{\P_0} \V_h)\big)$.
Since $\nabla_h q_h \in \P_0$ we deduce from proposition
\ref{prop:propadjh}
\begin{equation*}
\label{eq:infs5}
  -\big(q_h,\hbox{div}_h(\Pi_{\P_0} \V_h)\big)=(\Pi_{\P_0} \V_h,\nabla_h q_h)=(\V_h,\nabla_h q_h).
\end{equation*}
By splitting the last term we get
\begin{equation}
\label{eq:infs4}
  -\big(q_h,\hbox{div}_h(\Pi_{\P_0} \V_h)\big)=(\V,\nabla_h q_h)-(\V-\V_h,\nabla_h q_h).
\end{equation}
We bound the right-hand side of (\ref{eq:infs4}). Using
(\ref{eq:errintp1c}) and (\ref{eq:propv}) we have
\begin{equation*}
  |\V-\V_h|=|\V-\Pi_{\P^c_1} \V| \le C \, h \, \|\V\|_1 \le C \, h \,
  |q_h|.
\end{equation*}
Thus, using the Cauchy-Schwarz inequality, we get
\begin{equation*}
\label{eq:infs7}
  \left|(\V-\V_h,\nabla_h q_h)\right| \le C \, h \, |q_h|
   \, |\nabla_h q_h| \le C \, h \,  |q_h| \, \|q_h\|_{1,h}.
\end{equation*}
We estimate the other term as follows.
Integrating by parts we get
\begin{equation*}
(\V,\nabla_h q_h)=-(q_h,\hbox{div} \, \V)+ \sum_{K \in
{\mathcal{T}}_h} \sum_{\sigma \in {\mathcal{E}}_K} \int_{\sigma} q_h
\, (\V \cdot \N_{K,\sigma})
 \, d\sigma.
\end{equation*}
We have
  $-(q_h,\hbox{div} \, \V)=|q_h|^2$
thanks to  (\ref{eq:propv}). On the other hand
\begin{equation*}
\sum_{K \in {\mathcal{T}}_h} \sum_{\sigma \in {\mathcal{E}}_K}
\int_{\sigma} q_h \, (\V \cdot \N_{K,\sigma}) \, d\sigma
=\sum_{\sigma \in {\mathcal{E}}^{int}_h} \int_{\sigma} q_h \, (\V
\cdot \N_{K_\sigma,\sigma}) \, d\sigma
\end{equation*}
since $\V|_{\partial \Omega}=\mathbf{0}$. Using \cite{brenns} p.269
and (\ref{eq:propv}) we have
\begin{equation*}
 \left|\sum_{K \in {\mathcal{T}}_h} \sum_{\sigma \in {\mathcal{E}}_K} \int_{\sigma} q_h
 \, (\V \cdot \N_{K,\sigma}) \, d\sigma
 \right| \le C \, h \, \|\V\|_1 \, \|q_h\|_{1,h}
 \le C \, h \, |q_h|  \, \|q_h\|_{1,h}.
\end{equation*}
Hence we get
  $(\V,\nabla_h q_h) \ge (|q_h|-C \, h \, \|q_h\|_{1,h}) \, |q_h|$.
Thus we deduce from (\ref{eq:infs4})
\begin{equation}
\label{eq:infs8}
   -\big(q_h,\hbox{div}_h(\Pi_{\P_0} \V_h)\big) \ge (|q_h|-C \, h \, \|
   q_h\|_{1,h}) \, |q_h|.
\end{equation}
We now introduce the norm $\|.\|_h$. We have $\V_h=\Pi_{\P^c_1} \V
\in \P^c_1 \subset \H^1$. From \cite{eymgal} p. 776 we deduce
  $\|\Pi_{\P_0} \V_h\|_h \le C \, \|\V_h\|_1$.
Since $\Pi_{\P^c_1}$ is stable for the $\H^1$ norm, using
(\ref{eq:propv}), we get
\begin{equation*}
  \|\V_h\|_1=\|\Pi_{\P^c_1} \V\|_1 \le \|\V\|_1 \le C \, |q_h|.
\end{equation*}
Therefore
  $\|\Pi_{\P_0} \V_h\|_h \le C \, |q_h|$.
Using this  inequality in (\ref{eq:infs8}) we obtain that there
exists $C_1>0$ and $C_2>0$ such that
\begin{equation*}
   -\big(q_h,\hbox{div}_h(\Pi_{\P_0} \V_h)\big) \ge \left(C_1 \, |q_h|
   -C_2 \, h \, \|q_h\|_{1,h} \right) \,  \|\Pi_{\P_0} \V_h\|_h.
\end{equation*}
We deduce from this
\begin{equation*}
\sup_{\V_h \in \P_0\backslash \{\mathbf{0}\}}
    -\frac{(q_h,\hbox{\rm div}_h \, \V_h)}{\|\V_h\|_h} \ge C_1 \, |q_h|
    -C_2 \, h \, \|q_h\|_{1,h}.
\end{equation*}
Let us combine this result with lemma \ref{lemma:infsd}. Since
\begin{equation*}
  \forall \, t \ge 0 \, , \hspace{.5cm} \max \big(C \, t \, , \, C_1 \, |q_h| -C_2 \, t \big)
  \ge \frac{C \, C_1}{C+C_2} \, |q_h| \, ,
\end{equation*}
we finally get the result. \qed

\subsection{Properties of the discrete laplacian}
\label{subsec:proplapd}

We recall from \cite{zimm1} the coercivity of the laplacian
operator.
\begin{prop}
\label{prop:coerlap}
\noindent  For all $\U_h \in \P_0$ and $\V_h \in \P_0$ we have
  \begin{equation*}
  -(\Deltat_h \U_h,\U_h)=\|\U_h\|^2_h \, , \hspace{1cm}
    -(\Deltat_h \U_h,\V_h) \le \|\U_h\|_h \, \|\V_h\|_h.
  \end{equation*}
\end{prop}


\section{Stability of the scheme}
\label{sec:stab}

We first prove an estimate for the computed velocity (theorem
\ref{theo:estv}).
 We show  a similar result for the increments in time (lemma
 \ref{lem:estiv}). Using the inf-sup condition (proposition \ref{prop:condinfs}),
  we infer from it some estimates on the pressure (theorem \ref{theo:estp}).


\vspace{-.2cm}

\begin{lem}
\label{lemma:orthp} For all $m\in\{0,\dots,N\}$ et
$n\in\{0,\dots,N\}$ we have
\begin{equation*}
     (\U^m_h,\nabla_h p^n_h)=0  \, , \hspace{1cm}
     |\U^m_h|^2-|\Ut^m_h|^2+|\U^m_h-\Ut^m_h|^2=0.
\end{equation*}
\end{lem}
\noindent {\sc Proof.} First, using propositions \ref{prop:umrt0}
and \ref{prop:propadjh}, we get
 $(\U^m_h,\nabla_h p^n_h)=-(p^n_h,\hbox{div}_h \U^m_h)=0$.
Also, we deduce from   (\ref{eq:projib})
\begin{equation*}
 2 \, (\U^m_h,\U^m_h-\Ut^m_h)
=-\frac{4 \, k}{3} \big(\U^m_h,\nabla_h(p^m_h-p^{m-1}_h)\big)=0.
\end{equation*}
Using the  algebraic identity $2 \, a \, (a-b)=a^2-b^2+(a-b)^2$ we
get
\begin{equation*}
\hspace{4.3cm} 2 \,
(\U^m_h,\U^m_h-\Ut^m_h)=|\U^m_h|^2-|\Ut^m_h|^2+|\U^m_h-\Ut^m_h|^2=0.
\hspace{4.3cm} \qed
\end{equation*}

\vspace{.2cm}

\noindent We introduce the following hypothesis on the initial data.
\begin{equation*}
  {\bf (H1)}  \hspace{.2cm} \hbox{ There exists $C>0$ such that }
    \hspace{.2cm} |\U^0_h|+|\U^1_h|+k|\nabla_h p^1_h| \le C.
\end{equation*}
Hypothesis {\bf (H1)} is fulfilled if we set
$\U^0_h=\Pi_{\mathbf{RT_0}} \U_0$ and we use a semi-implicit Euler
 scheme to compute  $\U^1_h$. We have the following stability result.
\begin{theo}
\label{theo:estv} We assume  that the initial values  of the scheme
fulfill {\bf (H1)}.  For all $m \in \{2,\dots,N\}$ we have
\begin{equation}
\label{eq:estl2h1}
  |\U^m_h|^2+k\sum_{n=2}^m \|\Ut^n_h\|^2_h  \le C.
\end{equation}
\end{theo}
\noindent {\sc Proof.} Let $m \in\{2,\dots,N\}$
 and $n \in \{1,\dots,m-1\}$. Taking the scalar product of
 (\ref{eq:mombdf}) with
$4 \, k \, \Ut^{n+1}_h$ we get
\begin{align}
\label{eq:eqstabmultbdf}
    \left( \frac{3 \, \Ut^{n+1}_h-4 \, \U^n_h+\U^{n-1}_h}{2k}, 4 \, k \, \Ut^{n+1}_h \right)
&-\frac{4 \, k}{\hbox{Re}} \, (\Deltat_h \Ut^{n+1}_h, \Ut^{n+1}_h)
  \nonumber \\
 +4 \, k \, \b_h(2 \, \U^n_h -\U^{n-1}_h,\Ut^{n+1}_h,\Ut^{n+1}_h) &+4 \, k \, (\nabla_h p^n_h,\Ut^{n+1}_h)
=4 \, k \, (\F^{n+1}_h,\Ut^{n+1}_h).
\end{align}
First of all, using lemma \ref{lemma:orthp} and proceeding as in
\cite{guer}, we get
\begin{align*}
 4\, k \left(\Ut^{n+1}_h,\frac{3 \, \Ut^{n+1}_h-4 \,
\U^n_h+\U^{n-1}_h}{2 \, k}\right)
 &= |\U^{n+1}_h|^2-|\U^n_h|^2
    +|2 \, \U^{n+1}_h - \U^n_h|^2
  -|2\, \U^n_h - \U^{n-1}_h|^2  \\
  &+|\U^{n+1}_h -2 \, \U^n_h +\U^{n-1}_h|^2+6 \, |\Ut^{n+1}_h-\U^{n+1}_h|^2.
\end{align*}
According to  proposition \ref{prop:coerlap} we have
$-\frac{4 \, k}{\hbox{Re}} \, (\Deltat_h \Ut^{n+1}_h,
\Ut^{n+1}_h)=\frac{4 \, k}{\hbox{Re}} \, \|\Ut^{n+1}_h\|^2_h$.
Also, according to  lemma \ref{lemma:orthp} and (\ref{eq:projib})
\begin{eqnarray*}
\label{eq:estvbdf4}
4 \, k \, (\nabla_h p^n_h,\Ut^{n+1}_h)&=&4 \, k \, (\nabla_h p^n_h,\Ut^{n+1}_h-\U^{n+1}_h) \nonumber \\
&=&\frac{4 \, k^2}{3} \, ( |\nabla  p^{n+1}_h|^2 -
 |\nabla p^{n}_h|^2 - |\nabla  p^{n+1}_h -\nabla p^{n}_h |^2 ).
\end{eqnarray*}
Multiplying equation (\ref{eq:projib}) by $4 \, k \, \nabla_h
(p^{n+1}_h-p^n_h)$  and using the Young inequality we get
\begin{equation*}
\label{eq:estvbdf5} \frac{4 \, k^2}{3} \, |\nabla
(p^{n+1}_h-p^n_h)|^2 \le 3 \, |\U^{n+1}_h - \Ut^{n+1}_h|^2.
\end{equation*}
According to proposition \ref{prop:posbh}, we have
  $4 \, k \, \b_h(2 \, \U^n_h - \U^{n-1}_h,\Ut^{n+1}_h,\Ut^{n+1}_h) \ge 0$.
At last using the  Cauchy-Schwarz inequality, (\ref{eq:inpoinp0})
and (\ref{eq:propfh}) we have
\begin{equation*}
 4 \, k \, (\F^{n+1}_h, \Ut^{n+1}_h) \le 4 \, k \, |\F^{n+1}_h|\, |\Ut^{n+1}_h| \le C \, k \, \|\F\|_{{\mathcal{C}}(0,T;\L^2)}\, \|\Ut^{n+1}_h\|_h.
 \end{equation*}
Using the  Young inequality we get
\begin{equation*}
\label{eq:estvbdf3} 4 \, k \, (\F^{n+1}_h, \Ut^{n+1}_h) \le  3 \, k
\, \|\Ut^{n+1}_h\|^2_h + C \, k \,
\|\F\|^2_{{\mathcal{C}}(0,T;\L^2)}.
\end{equation*}
Thus we deduce from (\ref{eq:eqstabmultbdf})
\begin{eqnarray*}
\label{eq:estvbdf6}
 \hspace{-1cm} & & |\U^{n+1}_{h}|^2 - |\U^n_{h}|^2
+ |2 \, \U^{n+1}_h-\U^n_h|^2-|2 \, \U^n_h-\U^{n-1}_h|^2
 +|\U^{n+1}_h-2 \, \U^n_h+\U^{n-1}_h|^2 \nonumber \\
\hspace{-1cm} && +3 \, |\Ut^{n+1}_h - \U^{n+1}_h|^2
 +k \, \|\Ut^{n+1}_{h}\|^2_h +\frac{4 \, k^2}{3} \, (|\nabla_h p^{n+1}_h|^2 - |\nabla_h p^n_h|^2)  \le C \, k.
\end{eqnarray*}
Summing from $n=1$ to $m-1$ we have
\begin{eqnarray*}
\label{eq:estvbdf7} \hspace{-.6cm} && |\U^m_h|^2+|2 \,
\U^m_h-\U^{m-1}_h|^2 +3\sum_{n=1}^{m-1} |\Ut^{n+1}_h - \U^{n+1}_h|^2
+k\sum_{n=1}^{m-1} \|\Ut^{n+1}_h\|^2_h + \frac{4 \, k^2}{3} \, |\nabla_h p^m_h|^2  \nonumber \\
\hspace{-.6cm}&& \le C+4 \, |\U^1_h|^2+|2 \, \U^1_h - \U^0_h|^2+k^2
\, |\nabla_h p^1_h|^2.
\end{eqnarray*}
Using hypothesis {\bf (H1)} we get (\ref{eq:estl2h1}). \qed

\vspace{.2cm}

\noindent We now want to estimate the computed pressure.
From now on, we make the following hypothesis on the data
\begin{equation*}
    \F \in {\mathcal{C}}(0,T;\L^2) \, , \hspace{.6cm} \F_t \in L^2(0,T;\L^2) \, ,
  \hspace{.6cm} \U_0 \in \H^2 \cap \H^1_0 \, , \hspace{.6cm} \hbox{div} \, \U_0=0.
\end{equation*}
 One shows that if the data
$\U_0$ and $\F$ fulfill
 a compatibility condition \cite{heywood}
there exists a  solution  $(\U,p)$ to the  equations
(\ref{eq:mom})--(\ref{eq:incomp}) such that
\begin{equation*}
 \U \in {\mathcal{C}}(0,T;\H^2) \, , \hspace{.5cm} \U_t \in {\mathcal{C}}(0,T;\L^2) \, ,
  \hspace{.5cm} \nabla p\in {\mathcal{C}}(0,T;\L^2).
\end{equation*}
We introduce the following hypothesis on the initial values of the
scheme: there exists a constant $C>0$ such that
\begin{equation*}
  {\bf (H2)}    \hspace{.7cm}  |\U^0_h-\U_0| + \frac{1}{h} \, \|\U^1_h-\U(t_1)\|_{\infty}
      +|p^1_h-p(t_1)| \le C \, h \, , \hspace{.7cm} |\U^1_h-\U^0_h|\le C \, k.
\end{equation*}
One checks easily that this hypothesis implies {\bf (H1)}. We have
the following result.

\vspace{-.2cm}

\begin{lem}
\label{lem:estiv} We assume that the initial values of the scheme
fulfill {\bf (H2)}. Then there exists a constant $C>0$ such that for
all $m \in \{1,\dots,N\}$
\begin{equation*}
\label{eq:estincv}
  \frac{1}{k} \, |\U^m_h-\U^{m-1}_h| \le C .
\end{equation*}
\end{lem}
\noindent {\sc Proof.} Using proposition \ref{prop:consbh} one
proceeds as in \cite{zimm1}.  The difference lies in the way we
bound the term $\nabla_h p^1_h$. We use the splitting
\begin{equation*}
  p^1_h=(p^1_h-\Pi_{P^{nc}_1} p(t_1))+(\Pi_{P^{nc}_1}
  p(t_1)-p(t_1))+p(t_1).
\end{equation*}
Using an inverse inequality \cite{brenns}  we have
\begin{equation*}
\left|\nabla_h \left(p^1_h-\Pi_{P^{nc}_1} p(t_1)\right)\right| \le
\frac{C}{h} \, \left|p^1_h-\Pi_{P^{nc}_1} p(t_1)\right| \le
\frac{C}{h} \, \left(
\left|p^1_h-p(t_1)\right|+\left|p(t_1)-\Pi_{P^{nc}_1} p(t_1)\right|
\right).
\end{equation*}
Using (\ref{eq:estintp1nc}) and hypothesis {\bf (H2)} we get
\begin{equation*}
\left|\nabla_h \left(p^1_h-\Pi_{P^{nc}_1} p(t_1)\right)\right| \le C
\, \|p(t_1)\|_1 \le C \, \|p\|_{{\mathcal{C}}(0,T;\H^1)}.
\end{equation*}
According to (\ref{eq:estintp1nc}) we also have
  $\left| \nabla_h (p(t_1)-\Pi_{P^{nc}_1} p(t_1)) \right|
  \le C \, \|p(t_1)\|_1 \le C \, \|p\|_{{\mathcal{C}}(0,T;\H^1)}$.
Lastly
  $|\nabla p(t_1)| \le \|p\|_{{\mathcal{C}}(0,T;\H^1)}$.
Thus
  we get $|\nabla_h p^1_h| \le C$. \qed

\begin{theo}
\label{theo:estp} We assume that the initial values of the scheme
fulfull {\bf (H2)}. There exists a  constant $C>0$ such that for all
$m\in\{2,\dots,N\}$
\begin{equation*}
\label{eq:eqestp}
  k\sum_{n=2}^m |p^n_h|^2 \le C.
  \end{equation*}
\end{theo}
\noindent {\sc Proof.} Let $m \in \{2,\dots,N\}$. We set  $n=m-1$.
Using the  inf-sup condition (\ref{eq:infsd}) and proposition
\ref{prop:propadjh}, we get that there exists $\V_h \in
\P_0\backslash\{\mathbf{0}\}$ such that
\begin{equation}
\label{eq:stabpbdf456}
 C \, \|\V_h\|_h \, |p^{n+1}_h|
 \le -(p^{n+1}_h,\hbox{\rm div}_h \, \V_h)=(\nabla_h p^{n+1}_h,\V_h).
\end{equation}
Plugging (\ref{eq:projib}) into  (\ref{eq:mombdf}) we have
\begin{equation*}
\label{eq:stabpbdf455}
 \nabla_h p^{n+1}_h   = -\frac{3 \, \U^{n+1}_h - 4 \, \U^n_h+\U^{n-1}_h}{2 \, k}
+\frac{1}{\hbox{Re}} \, \Deltat_h \Ut^{n+1}_h -\bt_h(2 \, \U^n_h -
\U^{n-1}_h,\Ut^{n+1}_h)+ \F^{n+1}_h.
\end{equation*}
so that
\begin{align*}
 (\nabla_h p^{n+1}_h,\V_h)
 &= -\left(\frac{3 \, \U^{n+1}_h - 4 \, \U^n_h+\U^{n-1}_h}{2 \, k},\V_h\right)
+\frac{1}{\hbox{Re}}\left(\Deltat_h \Ut^{n+1}_h,\V_h\right) \\
 &-\b_h(2 \, \U^n_h-\U^{n-1}_h,\Ut^{n+1}_h,\V_h)
 + (\F^{n+1}_h, \V_h).
\end{align*}
Thanks to  proposition \ref{prop:stabbth} and theorem
\ref{theo:estv} we have
\begin{equation*}
 \left|\b_h(2 \, \U^n_h-\U^{n-1}_h,\Ut^{n+1}_h,\V_h)\right| \le \left(2 \, |\U^n_h|+|\U^{n-1}_h|\right)
  \, \|\Ut^{n+1}_h\|_h \, \|\V_h\|_h
  \le C \, \|\Ut^{n+1}_h\|_h \, \|\V_h\|_h.
\end{equation*}
According  to proposition \ref{prop:coerlap} we have
  $\left(\Deltat_h \Ut^{n+1}_h,\V_h\right) \le \|\Ut^{n+1}_h\|_h \, \|\V_h\|_h$.
Using the Cauchy-Schwarz inequality, (\ref{eq:inpoinp0}) and
(\ref{eq:propfh}) we have
\begin{equation*}
  (\F^{n+1}_h,\V_h) \le |\F^{n+1}_h| \, |\V_h| \le C \, |\V_h| \le C  \, \|\V_h\|_h
\end{equation*}
and in a similar way
\begin{equation*}
\left|\left(\frac{3 \, \U^{n+1}_h - 4 \, \U^n_h+\U^{n-1}_h}{2 \,
k},\V_h\right)\right| \le C \, \left|\frac{3 \, \U^{n+1}_h - 4 \,
\U^n_h+\U^{n-1}_h}{2 \, k}\right| \, \|\V_h\|_h.
\end{equation*}
Thus we get
\begin{equation*}
\label{eq:stabsbdf}
 (\nabla_h p^{n+1}_h,\V_h)    \le C +C\left( \frac{|3 \, \U^{n+1}_h - 4 \, \U^n_h+\U^{n-1}_h|}{2 \, k}
+\|\Ut^{n+1}_h\|_h\right) \|\V_h\|_h.
\end{equation*}
By comparing with (\ref{eq:stabpbdf456}) we get
\begin{equation*}
 |p^{n+1}_h|       \le C+C \left(
\frac{|3 \, \U^{n+1}_h - 4 \, \U^n_h+\U^{n-1}_h|}{2 \, k}
+\|\Ut^{n+1}_h\|_h  \right).
\end{equation*}
Squaring and  summing from $n=1$ to $m-1$ we obtain
\begin{equation*}
\label{eq:stabpbdf6}
  k\sum_{n=2}^{m} |p^{n}_h|^2
 \le   C+C \, k\sum_{n=1}^{m-1} \frac{|3 \, \U^{n+1}_h - 4 \, \U^n_h+\U^{n-1}_h|^2}{4 \, k^2}+
C \, k\sum_{n=1}^{m-1} \|\Ut^{n+1}_h\|^2_h.
\end{equation*}
The last term on the right-hand side is bounded, thanks to theorem
\ref{theo:estv}. And since
\begin{equation*}
3 \, \U^{n+1}_h - 4 \, \U^n_h+\U^{n-1}_h=3(\U^{n+1}_h - \U^n_h)
-(\U^n_h - \U^{n-1}_h) =3 \, \Del \U^{n+1}_h - \Del \U^n_h
\end{equation*}
we deduce from lemma \ref{lem:estiv}
\begin{equation*}
\label{eq:stabpbdf46}
 \hspace{4.3cm} k\sum_{n=1}^{m-1} \frac{|3 \, \U^{n+1}_h - 4 \, \U^n_h+\U^{n-1}_h|^2}{4 \, k^2}
 \le C \, k\sum_{n=1}^{m} \frac{|\Del \U^n_h|^2}{k^2} \le C.  \hspace{4.3cm}\qed
\end{equation*}




\begin{thebibliography}{2}
%





\bibitem{boicaya}  S. Boivin, F. Cayre, J. M Herard,
A finite volume method to solve the Navier-Stokes equations
 for incompressible flows on unstructured meshes,
\textit{Int. J. Therm. Sci.}  \textbf{39} (2000) 806--825.

\bibitem{brenns}
S. C. Brenner  and L. R. Scott, \textit{The mathematical theory of
finite element methods}, Springer, 2002.

\bibitem{brezzfor}
F. Brezzi and M. Fortin, \textit{Mixed and hybrid finite element
methods}, Springer-Verlag, 1991.

\bibitem{chorin}
J. Chorin, On the convergence of discrete approximations to the
Navier-Stokes equations, \textit{Math. Comp.} \textbf{23} (1969)
341--353.




\bibitem{eymgal}  R. Eymard, T. Gallou{\"e}t and R. Herbin,
\textit{Finite volume methods}, P.G. Ciarlet and J.L. Lions eds,
North-Holland, 2000.



\bibitem{herb3}  R. Eymard and R. Herbin, A staggered
finite volume scheme on general meshes for the Navier-Stokes
equations in two space dimensions, \textit{Int.J. Finite Volumes}
(2005).



\bibitem{eymard}
R. Eymard, J. C. Latch\'e and R. Herbin, Convergence analysis of a
colocated  finite volume scheme for the  incompressible
Navier-Stokes equations on  general 2 or 3D meshes, \textit{SIAM J.
Numer. Anal.} \textbf{45}(1) (2007) 1--36.




\bibitem{faure}
S. Faure, Stability of a colocated finite volume scheme for the
Navier-Stokes equations, \textit{Num. Meth. PDE} \textbf{21}(2)
(2005) 242--271.

\bibitem{giraultr} V. Girault and  P. A. Raviart,
\textit{Finite Element Methods for Navier-Stokes equations}: Theory
and Algorithms, Springer, 1986.

\bibitem{guer}
J. L. Guermond, Some implementations of projection methods for
Navier-Stokes equations, \textit{M2AN} \textbf{30}(5) (1996)
637--667.


\bibitem{heywood}
J. G. Heywood  and R. Rannacher, Finite element approximation of the
nonstationary Navier-Stokes problem. I. Regularity of solutions and
second-order error estimates for spatial discretization,
\textit{SIAM J. Numer. Anal.} \textbf{19}(26) (1982) 275--311.


\bibitem{kimchoi} D. Kim   and H. Choi, A
second-order time-accurate finite volume method for unsteady
incompressible flow on hybrid unstructured grids, \textit{J. Comp.
Phys.} \textbf{162}
 (2000) 411--428.




\bibitem{temam}
R. Temam, Sur l'approximation de la solution des \'equations de
Navier-Stokes par la m\'ethode de pas fractionnaires II,
\textit{Arch. Rat. Mech. Anal.} \textbf{33} (1969) 377--385.



\bibitem{zimm1}
S. Zimmermann, Stability of a colocated finite volume for the
incompressible Navier-Stokes equations, \textit{arXiv}:0704.0772
(2006).

\bibitem{zimm}
S. Zimmermann, \textit{\'Etude et impl\'ementation de m\'ethodes de
volumes finis pour les fluides incompressibles}, PhD, Blaise Pascal
University, France (2006).






















\end{thebibliography}
\end{document}